\documentclass[12pt,a4paper]{amsart}

\usepackage[latin1]{inputenc}
\usepackage{amsmath,amssymb,amsthm}

\numberwithin{equation}{section}

\newtheorem{theorem}[equation]{Theorem}

\newtheorem{lemma}[equation]{Lemma}
\newtheorem{corollary}[equation]{Corollary}

\theoremstyle{definition}
\newtheorem{definition}[equation]{Definition}
\newtheorem{example}[equation]{Example}
\newtheorem{remark}[equation]{Remark}

\newcommand{\N}{\mathbb N}

\newcommand{\R}{\mathbb R}
\newcommand{\Z}{\mathbb Z}
\newcommand{\C}{\mathbb C}

\newcommand{\sgn}{\operatorname{sgn}}

\newcommand{\bddlin}{\mathcal L}

\newcommand{\tr}{{\rm tr}}
\newcommand{\lin}{{\rm lin}}

\newcommand{\<}{\langle}
\renewcommand{\>}{\rangle}
\newcommand{\ket}[1]{\left|#1\right\rangle}
\newcommand{\bra}[1]{\left\langle#1\right|}
\newcommand{\kb}[2]{\left|#1\right\rangle\left\langle#2\right|}
\newcommand{\norm}[2]{\|#1\|_{#2}}

\title[Nonpositive sesquilinear form measures]{Diagonalization and representation results for nonpositive
sesquilinear form measures}
\author[Hyt\"onen \and Pellonp\"a\"a \and Ylinen]{}
\date{\today}

\begin{document}

\maketitle

\begin{center}

Tuomas Hyt\"onen\\
Department of Mathematics and Statistics\\ University of Helsinki\\ Gustaf H\"allstr\"omin katu 2b\\ FI-00014 Helsinki, Finland\\ \texttt{tuomas.hytonen@helsinki.fi}\\ \mbox{}\\

Juha-Pekka Pellonp\"a\"a\footnote{Corresponding author. Telephone: +358-2-333 5737, telefax: +358-2-333 5070.}\\

Department of Physics\\ University of Turku\\ FI-20014 Turku, Finland\\  \texttt{juhpello@utu.fi}\\ \mbox{}\\

Kari Ylinen\\
Department of Mathematics\\ University of Turku\\ FI-20014 Turku, Finland\\ \texttt{ylinen@utu.fi}

\end{center}


\subsection*{Abstract}
We study decompositions of operator measures and more general sesquilinear form measures $E$ into linear combinations of positive parts, and their diagonal vector expansions. The underlying philosophy is to represent $E$ as a trace class valued measure of bounded variation on a new Hilbert space related to $E$. The choice of the auxiliary Hilbert space fixes a unique decomposition with certain properties, but this choice itself is not canonical. We present relations to Naimark type dilations and direct integrals.

\subsection*{Mathematics Subject Classification}
47A70 (Primary); 47A07 (Secondary).

\subsection*{Keywords and phrases}
Sesquilinear form, operator measure, bounded variation, diagonalization, Naimark dilation, direct integral.




\section{Introduction}

The idea of rigged Hilbert spaces arises in attempts to develop
mathematically rigorous interpretations of the intuitively
appealing Dirac formalism of Quantum Mechanics. With the help of
generalized eigenvectors lying outside the Hilbertian state space,
one is able to write eigenvalue expansions, with formal similarity
to the finite-dimensional case, even for self-adjoint operators
with a continuous spectrum. By the Spectral Theorem, self-adjoint
operators may be identified with spectral
measures on the real line, and they are the mathematical
representatives of physical observables in the traditional von
Neumann approach to Quantum Mechanics. It is, however, well known
that this point of view becomes too restrictive already when
considering such basic physical examples as phase-like quantities 
(see e.g.\ \cite{JP}),
but they can still be incorporated into the mathematical formalism
by allowing more general positive operator measures in place of
spectral measures. Also for them, and for the yet larger class of
positive sesquilinear form measures, generalized eigenvalue
expansions have been obtained in the literature. See e.g.
\cite{HPY} and the references therein.

The present note is concerned with similar results for
sesquilinear form measures without any positivity conditions.
Besides purely mathematical interest, motivation comes from
important physical questions. Let us consider an example.

\begin{example}
Let $H$ be a complex Hilbert space with an orthonormal basis
$(e_n)_{n=0}^\infty$. Let $z\in\C\setminus\{0\}$ and define a {\it
coherent state}
$$
\psi_z:=e^{-|z|^2/2}\sum_{n=0}^\infty\frac{z^n}{\sqrt{n!}}e_n.
$$
It describes quasimonochromatic laser light (in a single-mode
quantum optical system), where $|z|$ is the energy parameter and $z/|z|\in\mathbb T$ is the phase parameter of
the laser light (see e.g.\ \cite{Per}). The vector $e_n$, the so-called {\it number state} or {\it Fock state},
describes an optical field which contains $n$ photons of the same
frequency.

A measurement of the phase parameter can be described by using a
phase shift covariant semispectral measure \cite[p.\ 23]{JP}
$$
E(X):=\sum_{m,n=0}^\infty c_{mn}\int_X w^{n-m}d\mu(w)\kb{e_m}{e_n}
$$
where $X$ is a Borel set of $\mathbb T$, $\mu$ is the normalized
Haar measure of $\mathbb T$ and $(c_{mn})$ is a positive
semidefinite complex matrix with the unit diagonal; the
probability of getting a value $w$ from a set $X$ when the system
is prepared in a state $\psi_z$ is thus
$\langle\psi_z|E(X)\psi_z\rangle$.

In realistic physical situations we cannot produce arbitrarily
high photon numbers, that is, we cannot prepare number states
$e_n$ for an arbitrarily large $n$. In fact, as of 2004, a method described in \cite{VBW} ``still remains the only experiment in principle capable of providing an arbitrary Fock state (at least up to $n=4$) on demand''.
But still we need the whole
Hilbert space $H$ to define coherent states. Hence, we can relax
the definition of $(c_{mn})$: we need only assume that the first,
say, $10\times 10$ block of $(c_{mn})$ is positive semidefinite
(so that we get probability distributions also for superpositions of number states $e_n$, $n\le 10$).
Moreover, we assume that for any coherent state $\psi_z$ with sufficiently low energy $|z|\le r\in\R$ we can
define a probability measure 
$X\mapsto e^{-|z|^2}\sum_{m,n=0}^\infty c_{mn}\int_X w^{n-m}d\mu(w)\frac{z^{n-m}}{\sqrt{n!m!}}$.
Further restrictions can be imposed, if we assume that some superpositions of
coherent states can be measured.
If $(c_{mn})$ is not assumed to be positive semidefinite,
then $E(X)$ may be a nonpositive operator, or even a sesquilinear form on
$V=\lin\{e_n\,|\,n=0,1,\ldots\}$, for some $X$. Then the mapping $X\mapsto E(X)$ can be understood as a (nonpositive) {\it sesquilinear form measure}.
It can be shown that some (phase shift covariant) sesquilinear form measures give more accurate phase distributions in coherent states than any (covariant) positive semispectral measures \cite{JP:IV}. 
\end{example}

The sesquilinear form measures we study here generalize operator
measures which have already received a fair amount of attention in
the mathematical literature. For example, we may quote a
well-known decomposition result from \cite{Pa02}, pp. 104--105: A
regular Borel operator measure on a compact Hausdorff space (with
values in the space of bounded operators on a Hilbert space) is,
as a consequence of Wittstock's decomposition theorem,  completely
bounded if, and only if, it can be expressed as a linear
combination of positive operator valued measures. In this paper an
analogous decomposition problem in the setting of sesquilinear
form measures is in a central role. We consider a $\sigma$-algebra
$\Sigma$, a vector space $V$ with a countable Hamel basis, and
measures $E:\Sigma\to S(V)$ where $S(V)$ is the space of
sesquilinear forms on $V$. This generalizes the more standard
setting of operator measures in the context of a separable Hilbert
space, and it turns out that our more flexible framework yields
new information even there: An operator measure may be decomposed
into a linear combination of positive parts without the condition
of complete boundedness. Of course there is a price to pay: these
positive parts are not necessarily operator valued but only
sesquilinear form valued. While this on the one hand may be seen
as a drawback, on the other hand it highlights the usefulness of
general sesquilinear form measures.

The paper is organized as follows. After the setting is explained
in Section 2, the next section establishes a connection with
operator measures taking their values in the trace class
${\mathcal L}^1(H)$ of a separable Hilbert space. Since ${\mathcal
L}^1(H)$ has the Radon-Nikod\'ym property, a sesquilinear form
measure can be expressed in terms of integrating an ${\mathcal
L}^1(H)$-valued density function with respect to a basic positive
scalar measure. In Section 4 the desired decomposition is effected
by utilizing the operator density found in Section 3. While the
basic idea is straightforward enough, one must take care of rather
delicate measurability issues. To this end, a classical result of
Kuratowski and Ryll-Nardzewski on measurable selectors is used.
The final Section 5 deals with an analogue of the Naimark dilation
theorem: The decomposition of a sesquilinear form measure into
positive parts also yields a spectral dilation in a generalized
sense involving a unitary operator $W$ on the dilation space where
the spectral measure acts. The characateristic feature of $W$ is
that $W^4=I$. The paper concludes with a remark on formulating
the dilation result in terms of a direct integral representation.


\section{Basics}

We write $\Z_+:=\{1,2,3,\ldots\}$, $\N:=\{0\}\cup\Z_+$ and
$\Z_-:=\Z\setminus\N$. For $p>0$ and $\mathcal I$ an index set,
$\ell^p(\mathcal I)$ is the space of the complex families
$c=(c_n)_{n\in\mathcal I}$ such that $\sum_{n\in\mathcal
I}|c_n|^p<\infty$. 

Let $V$ be a vector space. The scalar field is always $\C$. A
mapping $\Phi:V\times V\to\C$ is called a {\em sesquilinear form}
(SF), if it is antilinear (i.e., conjugate linear) in the first
and linear in the second variable.
 It is {\it symmetric} if $\Phi(\phi,\psi)=\overline{\Phi(\psi,\phi)}=:\Phi^*(\phi,\psi)$ and
 {\em positive} if $\Phi(\phi,\phi)\geq0$ for all $\phi,\,\psi\in V$.
Any positive SF is symmetric, and any SF $\Phi$ is a linear
combination of two symmetric SFs:
\begin{equation}\label{eq:saParts}
  \Phi=\frac{1}{2}(\Phi+\Phi^*)+\frac{i}{2}(i\Phi^*-i\Phi).
\end{equation}
We let $S(V)$ (resp.\ $PS(V)$) denote the set of sesquilinear
forms (resp.\ positive sesquilinear forms) on $V\times V$.

Our basic reference on measure and (vector) integration is
\cite{DS}. Measurability means $\mu$-measurability where $\mu$ is
a fixed positive measure. Let $(\Omega,\Sigma)$ be a measurable
space, i.e., $\Sigma$ is a $\sigma$-algebra of subsets of
$\Omega$.

\begin{definition}\label{def:PSFM}
Let $E:\Sigma\to S(V)$ be a mapping and denote $E(X)=E_X$ for
$X\in \Sigma$. We call $E$ a {\em sesquilinear form measure (SFM)}
 if  the mapping $X\mapsto
E_X(\phi,\psi)$ is $\sigma$-additive, i.e. a complex measure, for
all $\phi,\,\psi\in V$. If in addition $E(X)$ is symmetric (resp.\
positive) for all $X\in\Sigma$, $E$ is called a {\em symmetric }
(resp.\ {\em positive) sesquilinear form measure}.
\end{definition}

The inner product of any Hilbert space $H$ is linear in the second
variable and denoted by $\<\,\cdot\,|\,\cdot\,\>$. We let
$\bddlin(H)$ stand for the bounded linear operators on $H$,
$\bddlin_s(H)\subset\bddlin(H)$ for the self-adjoint operators,
and $\bddlin_+(H)\subset\bddlin_s(H)$ for the positive ones. The
trace class is denoted by $\bddlin^1(H)$, and
$\bddlin^1_s(H):=\bddlin^1(H)\cap\bddlin_s(H)$,
$\bddlin^1_+(H):=\bddlin^1(H)\cap\bddlin_+(H)$.

\begin{definition}
Let $H$ be a Hilbert space and $E_0:\Sigma\to\bddlin(H)$ a
mapping. We call $E_0$ an \emph{operator measure (OM)} if it is
weakly $\sigma$-additive, i.e.\ the mapping
$X\mapsto\<\phi|E_0(X)\psi\>$ is $\sigma$-additive for all
$\phi,\,\psi\in H$. If in addition
$E_0(\Sigma)\subset\bddlin_s(H)$ (resp.\
$E_0(\Sigma)\subset\bddlin_+(H)$) we say that $E_0$ is a {\em
self-adjoint} (resp.\ {\em positive) operator measure}, and if
$E_0(X)^2=E_0(X)=E_0(X)^*$ for all $X\in \Sigma$, $E_0$ is called
a {\em projection measure}. An OM $E_0:\Sigma\to \bddlin(H)$ is
called \emph{normalized} if $E_0(\Omega)=I$, the identity operator
on $H$. A normalized positive OM is also called a
\emph{semispectral measure} and a normalized projection measure a
\emph{spectral measure}.
\end{definition}

Every (self-adjoint or positive) OM $E_0$ can be identified with a
(symmetric or positive) SFM $E$ by setting
$E_X(\phi,\psi):=\<\phi|E_0(X)\psi\>$.

\section{Reduction to trace-class operator measures}

For the rest of the note, we assume that $V$ has a countably
infinite Hamel basis $(e_n)_{n=0}^{\infty}$, and $H$ is the
Hilbert space completion of $V$ such that
 $(e_n)_{n=0}^{\infty}$ is an orthonormal basis of $H$.
For any SF $\Phi$ on $V$ we write (formally)
$$
\Phi=\sum_{m,n=0}^\infty \Phi_{mn}\kb{e_m}{e_n}
$$
where $\Phi_{mn}:=\Phi(e_m,e_n)$. If $\Phi$ is bounded with
respect to the norm of $H$, it determines a unique bounded linear
operator $\tilde\Phi\in\bddlin(H)$ satisfying
$\<\phi|\tilde\Phi\psi\>=\Phi(\phi,\psi)$. Then the series above
is not just formal; when $\kb{e_m}{e_n}$ denotes as usual the rank
one operator $\phi\mapsto\<e_n|\phi\>e_m$,  the series converges
with respect to the weak operator topology to $\tilde\Phi$. We may
identify $\Phi$ and $\tilde\Phi$, and then $\Phi_{mn}=\<e_m|\Phi
e_n\>$.

\begin{lemma}\label{lem:tracetest}
Let $\Phi\in S(V)$ be represented by an infinite matrix
$(\Phi_{mn})_{m,n=0}^{\infty}\in\ell^1(\N\times\N)$. Then $\Phi$
has a unique extension $\Phi\in\bddlin^1(H)$ and
$\|\Phi\|_{\bddlin^1(H)}\leq\sum_{m,n=0}^{\infty}|\Phi_{mn}|$.
\end{lemma}

\begin{proof}
Since $\bddlin^1(H)$ is the dual of the space of finite rank
operators on $H$, the first claim is equivalent to requiring that
$
  \sup|\tr(\Phi \Lambda)|<\infty
$ where $\Phi$ is interpreted as a matrix $(\Phi_{mn})$ and
$\Lambda$ ranges over the matrices $(\Lambda_{mn})$ of finite rank
operators of norm $\le1$. But
\[
  |\tr(\Phi\Lambda)|=\Big|\sum_{m,n=0 }^\infty \Phi_{mn}\Lambda_{nm}\Big|
  \leq\sum_{m,n=0}^\infty |\Phi_{mn}|\cdot|\Lambda_{nm}|
  \leq\sum_{m,n=0}^\infty |\Phi_{mn}|,
\]
since $|\Lambda_{nm}|\leq\|\Lambda\|_{\bddlin(H)}\leq 1$.
\end{proof}

\begin{theorem}\label{thm:traceRep}
For any SFM $E:\Sigma\to S(V)$ there exist an
$\bddlin^1(H)$-valued measure $F$ of bounded variation, and an
injective operator $D\in\bddlin_+(H)$ such that $DV=V$ and
\begin{equation*}
  E_X(D\phi,D\psi)=\<\phi|F(X)\psi\>,\qquad\phi,\ \psi\in V.
\end{equation*}
There further exist a finite positive measure
\(\mu:\Sigma\to[0,\infty)\) and a function
\begin{equation*}
  T\in L^1(\Omega,\Sigma,\mu;\bddlin^1(H))
\end{equation*}
such that
\begin{equation*}
  E_X(D\phi,D\psi)=\int_X \<\phi|T(\omega)\psi\>d\mu(\omega),\qquad\phi,\ \psi\in V.
\end{equation*}
\end{theorem}

Defining 
$
  C_{\omega}(\phi,\psi):=\<D^{-1}\phi|T(\omega)D^{-1}\psi\>,
$
we also obtain the integral representation
\begin{equation*}
  E_X(\phi,\psi)=\int_X C_{\omega}(\phi,\psi)d\mu(\omega),\qquad\phi,\ \psi\in V.
\end{equation*}

\begin{proof}
We denote $E_{mn}(X):=E_X(e_m,e_n)$, and write $|E_{mn}|(X)$ for
its total variation on $X$.

Choose any bounded positive sequence \((d_m)_{m=0}^{\infty}\) such
that
\begin{equation*}
  \delta:=\sum_{m,n=0}^{\infty}d_m d_n|E_{mn}|(\Omega)<\infty.
\end{equation*}
For example, we may take
$d_m=\alpha_m/\max\{1,\sqrt{|E_{kl}|(\Omega)}\,|\,0\le k,l\le m\}$
where $(\alpha_m)_{m=0}^{\infty}$ is any summable positive
sequence.

Let $D$ be the diagonal operator
\begin{equation*}
  D\phi:=\sum_{n=0}^{\infty}d_n\ket{e_n}\<e_n|\phi\>.
\end{equation*}
Then for \(\phi,\ \psi\in V\),
\begin{equation}\label{eq:defF}
  E_X(D\phi,D\psi)
  =\sum_{m,n=0}^{\infty}\<\phi|e_m\>d_m d_n E_{mn}(X)\<e_n|\psi\>
  =:\<\phi|F(X)\psi\>,
\end{equation}
and we have \(F(X)\in\bddlin^1(H)\) with
\(\|F(X)\|_{\bddlin^1(H)}\leq\delta\) by
Lemma~\ref{lem:tracetest}. If \((X_k)_{k=0}^{\infty}\) is any
countable partition of \(X\subset\Omega\), then
\begin{equation}\label{eq:bddVarEst}
  \sum_{k=0}^{\infty}\|F(X_k)\|_{\bddlin^1(H)}
  \leq\sum_{k=0}^{\infty}\sum_{m,n=0}^{\infty}d_m d_n|E_{mn}(X_k)|
  \leq\sum_{m,n=0}^{\infty}d_m d_n|E_{mn}|(X)\leq\delta.
\end{equation}
This justifies the computation
\begin{equation*}
  \sum_{k=0}^{\infty}F(X_k)
  =\sum_{k=0}^{\infty}\sum_{m,n=0}^{\infty}d_m d_n E_{mn}(X_k)\ket{e_m}\bra{e_n}
  =\sum_{m,n=0}^{\infty}d_m d_n E_{mn}(X)\ket{e_m}\bra{e_n}
  =F(X),
\end{equation*}
which shows that $F$ is $\sigma$-additive, and
\eqref{eq:bddVarEst} with $X=\Omega$ also shows that $F$ is of
bounded variation.

For the measure $\mu$ one can take any finite positive measure
with respect to which the vector measure $F$, or equivalently $E$, is absolutely continuous (i.e., whenever $\mu(X)=0$, we have also $F(X)=0$, or equivalently $E_X=0$ as a sesquilinear form). To be specific, we take $\mu$ to be the total variation of $F$,
\[
  |F|(X):=\sup\sum_{k=1}^N\|F(X_k)\|_{\bddlin^1(H)}
\]
where the supremum is over all finite $\Sigma$-partitions of $X$. 
As in the proof of Proposition 7.1 of \cite{HPY},
the existence of
$T$ then follows from the vector-valued Radon--Nikod\'ym theorem,
since $\bddlin^1(H)$ (as a separable dual space) has the
Radon--Nikod\'ym property.
\end{proof}

\begin{remark}\label{rem:HD}
The above theorem shows that a sesquilinear form measure on $V$ can always be viewed as an operator measure on a new Hilbert space. In fact, let us denote by $H_D$ the range of $D\in\bddlin(H)$ equipped with the inner product $\<\eta|\theta\>_D:=\<D^{-1}\eta|D^{-1}\theta\>$ and the induced norm. Then $D:H\to H_D$ is an isometric Hilbert space isomorphism. Observe that in \eqref{eq:defF} the series in the middle is absolutely convergent, and the right-hand side makes sense, for all $\phi,\,\psi\in H$. Thus $E_X$ extends continuously to a sesquilinear form on $H_D$, and for $\eta=D\phi,\,\theta=D\psi\in H_D$ we have
\begin{equation*}\begin{split}
  E_X(\eta,\theta) &=\<D^{-1}\eta|F(X)D^{-1}\theta\>
   =\<D^{-1}\eta|D^{-1}DF(X)D^{-1}\theta\> \\
  &=\<\eta|DF(X)D^{-1}\theta\>_D=:\<\eta|\tilde{E}(X)\theta\>_D.
\end{split}\end{equation*}
Due to the operator-ideal property of the trace class, we find that $X\in\Sigma\mapsto \tilde{E}(X)=DF(X)D^{-1}$ is an $\bddlin^1(H_D)$-valued measure of bounded variation. By the Radon--Nikod\'ym theorem, it can be written as
\[
  \tilde{E}(X)=\int_X S(\omega)d\mu(\omega),\qquad
  S\in L^1(\Omega,\Sigma,\mu;\bddlin^1(H_D)).
\]
\end{remark}

\begin{remark} In the rest of the paper we take $\mu$ to be
the measure constructed in the above proof. Assume now that
$\{\omega\}\in\Sigma$ for all $\omega\in\Omega$. If we let
$\mu=\mu_1+\mu_2$ be the decomposition of $\mu$ as the sum of a
discrete measure $\mu_1$ and a continuous measure $\mu_2$, the
integral formula in the above theorem may be used to decompose
$E$ as $E=E_1+E_2$ where $E_1$ is a discrete SFM, i.e.\ vanishes
outside a countable set, and the SFM $E_2$ is continuous, i.e.,
vanishes at every singleton. Clearly such a decomposition is
unique. \end{remark}

\section{Diagonalization; positive and negative parts}\label{sec:diag}

By formula~\eqref{eq:saParts} we may decompose the measures $E$
and $F$ as well as the operator density $T$ into linear
combinations of two symmetric parts, and by linearity the
representation formulae of Theorem~\ref{thm:traceRep} remain
true for these parts. In this section we obtain a further
decomposition of these symmetric parts. We will need the following
classical result on \emph{measurable selectors} from
\cite{Kuratowski}; it is also stated in \cite{DPZ}, Lemma~1.9:

\begin{lemma}\label{lem:Kuratowski}
Let $E$ be a compact metric space and let
$\psi:E\times\Omega\to\R$ be a mapping such that $\psi(x,\cdot)$
is measurable for arbitrary $x\in E$ and $\psi(\cdot,\omega)$ is
continuous for arbitrary $\omega\in\Omega$. Then there exists a
measurable $\xi:\Omega\to E$ such that
\[
  \psi(\xi(\omega),\omega)=\max_{x\in E}\psi(x,\omega),\hspace{0.3cm}
  \omega\in\Omega.
\]
\end{lemma}

\begin{corollary}\label{cor:Kuratowski}
Let $T:\Omega\to\bddlin^1(H)$ be a measurable function. Then there
exists a measurable $\Phi:\Omega\to \bar{B}_H$, the closed unit
ball of $H$, such that
\[
  |\<\Phi(\omega)|T(\omega)\Phi(\omega)\>|
  =\max_{\phi\in\bar{B}_H}|\<\phi|T(\omega)\phi\>|,\hspace{0.3cm}
  \omega\in\Omega.
\]
\end{corollary}

\begin{proof}
It is well known that the unit ball $\bar{B}_H$ of a separable
Hilbert space, when equipped with the weak topology, is a compact
metrizable space. We consider the mapping
\[
  \psi:\bar{B}_H\times\Omega\to\R,
  (\phi,\omega)\mapsto\<\phi|T(\omega)\phi\>,
\]
and it suffices to check the conditions of
Lemma~\ref{lem:Kuratowski}.

That $\psi(\phi,\cdot)$ is measurable is clear from the
assumptions. To see that $\psi(\cdot,\omega)$ is continuous,
denote $\Lambda:=T(\omega)\in\bddlin^1(H)$. Assume first that
$\Lambda=\ket{\psi_1}\bra{\psi_2}$ has rank $1$. The mappings
$\phi\mapsto\<\psi_i|\phi\>$ are obviously continuous in the
topology in question, and so is their product. In general, we have
\[
  \Lambda=\sum_{k=1}^{\infty}\ket{\psi_k}\bra{\rho_k},
  \qquad\sum_{k=1}^{\infty}\|\psi_k\|\cdot\|\rho_k\|<\infty.
\]
Since uniformly convergent series of continuous functions are
continuous, we have reached the conclusion.
\end{proof}

We can now prove a measurable diagonalization of an
$\bddlin^1_s(H)$-valued function. The proof follows closely the
same pattern as the special case for $\bddlin^1_+(H)$-valued
functions given in \cite{DPZ}, Proposition~1.8, but we include the
details for the reader's convenience.

\begin{theorem}\label{thm:diagonal}
Given a measurable function $T:\Omega\to\bddlin_s^1(H)$, there
exist measurable functions $\phi_k:\Omega\to H$ and
$\lambda_k:\Omega\to\R$, $k\in\Z_+$, such that for any fixed
$\omega\in\Omega$ there holds
\[\begin{split}
  &\<\phi_k(\omega)|\phi_{\ell}(\omega)\>=\delta_{k\ell},\qquad
  |\lambda_k(\omega)|\geq|\lambda_{\ell}(\omega)|\quad\text{if}\quad k\leq\ell,\\
  &T(\omega)=\sum_{k=1}^{\infty}\lambda_k(\omega)
  \ket{\phi_k(\omega)}\bra{\phi_k(\omega)},\qquad
  \|T(\omega)\|_{\bddlin^1(H)}
  =\sum_{k=1}^{\infty}\lambda_k(\omega).
\end{split}\]
\end{theorem}

\begin{proof}
This representation of $T(\omega)$ for each fixed
$\omega\in\Omega$ is just the usual spectral representation, but
the point is to obtain this with a measurable dependence on
$\omega$. To see this, we recall an algorithm for computing the
spectral representation. An eigenvalue $\lambda$ of
$\Lambda\in\bddlin^1_s(H)$ of largest modulus satisfies
\[
  |\lambda|=\max_{\phi\in\bar{B}_H}|\<\phi|\Lambda\phi\>|,
\]
and any $\phi\in\bar{B}_H$, which gives the maximum, is an
eigenvector related to $\pm\lambda$. By
Corollary~\ref{cor:Kuratowski}, there is a measurable function
$\phi_1:\Omega\to\bar{B}_H$ such that
\[
  \lambda_1(\omega):=\<\phi_1(\omega)|T(\omega)\phi_1(\omega)\>,
\]
which is also a measurable function of $\omega$ by the above
formula, is an eigenvalue of $T(\omega)$ of maximal modulus, with
the eigenvector $\phi_1(\omega)$.

We then repeat the same procedure with
$T_1(\omega):=T(\omega)-\lambda_1(\omega)\ket{\phi_1(\omega)}\bra{\phi_1(\omega)}$
in place of $T(\omega)$, obtaining new measurable functions
$\lambda_2(\omega)$ and $\phi_2(\omega)$. Proceeding inductively,
we obtain sequences of measureable functions
$(\lambda_k(\omega))_{k=1}^{\infty}$ and
$(\phi_k(\omega))_{k=1}^{\infty}$. At each fixed
$\omega\in\Omega$, these give the spectral decomposition of
$T(\omega)$ by standard results about compact selfadjoint
operators.
\end{proof}

It is now also easy to separate the positive and negative parts of
the operator density in a measurable way:

\begin{corollary}\label{cor:diag}
Given a measurable function $T:\Omega\to\bddlin_s^1(H)$, there
exist measurable functions $g_k:\Omega\to H$,
$k\in\Z\setminus\{0\}$, such that for any fixed $\omega\in\Omega$
there holds
\[\begin{split}
  &\<g_k(\omega)|g_{\ell}(\omega)\>
    =\delta_{k\ell}\|g_k(\omega)\|^2,\\
  &\|g_k(\omega)\|\geq\|g_{\ell}(\omega)\|\quad\text{if}\quad 0<k<\ell\quad\text{or}\quad 0>k>\ell,\\
  &T(\omega)=\sum_{k\in\Z\setminus\{0\}}\sgn(k)
  \ket{g_k(\omega)}\bra{g_k(\omega)},\qquad
  \|T(\omega)\|_{\bddlin^1(H)}
  =\sum_{k\in\Z\setminus\{0\}}\|g_k(\omega)\|^2.
\end{split}\]
\end{corollary}

\begin{proof}
With the notation of Theorem~\ref{thm:diagonal}, we define
the measurable functions
\[\begin{split}
  n_0(\omega) &:=0,\qquad{g_0(\omega)}:=0,\\
  n_{\pm k}(\omega) &:=\inf\{n\in\Z_+| n>n_{\pm(k-1)}(\omega),\
     \pm\lambda_n(\omega)>0\},\qquad k\in\Z_+ \\
  g_{\pm k}(\omega) &:=
     |\lambda_{n_{\pm k}(\omega)}(\omega)|^{1/2}
     {\phi_{n_{\pm k}(\omega)}(\omega)},\qquad k\in \Z_+,
\end{split}\]
where it is understood that $\inf\emptyset:=\infty$ and
$\lambda_{\infty}(\omega):=0=:{\phi_{\infty}(\omega)}$.
\end{proof}

\begin{corollary}\label{cor:Tpm}
Given a measurable function $T:\Omega\to\bddlin_s^1(H)$, there
exists a pair of measurable functions
$T^{\pm}:\Omega\to\bddlin_+^1(H)$, such that for any fixed
$\omega\in\Omega$ we have
\begin{itemize}
\item[(i)]
$
 T(\omega)=T^+(\omega)-T^-(\omega),
$
\item[(ii)]
$ T^+(\omega)T^-(\omega)=0, $ and
\item[(iii)]
$
\norm{T(\omega)}{\bddlin^1}=\norm{T^+(\omega)}{\bddlin^1}+\norm{T^-(\omega)}{\bddlin^1}.
$
\end{itemize}
Moreover, if (i) and (ii), or alternatively (i) and (iii), hold
for all $\omega\in\Omega$, the functions $T^+$ and $T^-$ are
uniquely determined.

\end{corollary}

\begin{proof}
For existence, it suffices to set
\begin{equation*}
  T^{\pm}(\omega):=\sum_{k\in\Z_{\pm}}\ket{g_k(\omega)}\bra{g_k(\omega)}.
\end{equation*}
The  uniqueness statement assuming (i) and (ii) follows e.g. from
Corollary 2.10 in \cite{SZ79}. Assuming (i) and (iii), the
uniqueness claim is a consequence of Theorem 4.2 in \cite{Tak79},
since ${\mathcal L}^1(H)$ with its norm and order may be
identified with the predual of $\mathcal L(H)$.
\end{proof}

In the case of a {\em symmetric} SFM $E$, its trace-class density
$T$ is self-adjoint operator valued and, using the above
corollaries, we get
\begin{equation}\label{eq:dk}\begin{split}
  E_X(\phi,\psi)
  &=\int_X\<D^{-1}\phi\,|\,[T^+(\omega)-T^-(\omega)]D^{-1}\psi\>d\mu(\omega)\\
  &=\int_X\sum_{k\in\Z\setminus\{0\}}\sgn(k)\<D^{-1}\phi|g_k(\omega)\>
     \<g_k(\omega)|D^{-1}\psi\>d\mu(\omega) \\
  &=\int_X\sum_{k\in\Z\setminus\{0\}}\sgn(k)\<\phi|d_k(\omega)\>
     \<d_k(\omega)|\psi\>d\mu(\omega),
\end{split}\end{equation}
where we have defined
\[
  {d_k(\omega)}:=D^{-1}{g_k(\omega)}\in H_{D^{-1}},
\]
and $H_{D^{-1}}$ is the Hilbert space consisting of all the formal sums $\sum_{n=0}^{\infty}c_n e_n$ such that $\sum_{n=0}^{\infty}d_n^2|c_n|^2<\infty$. Note that we have a Hilbert space triplet $H_D\subset H\subset H_{D^{-1}}$, where $H_{D^{-1}}$ is the topological antidual of $H_D$. Note that the conclusion of \eqref{eq:dk} could also have been reached by applying Corollary~\ref{cor:diag} to (the symmetric parts of) the function $S:\Omega\to L^1(\Omega,\Sigma,\mu;\bddlin^1(H_D))$ from Remark~\ref{rem:HD}.

Denoting
\[
  E^{\pm}_X(\phi,\psi):=
   \int_X\sum_{k\in\Z_{\pm}}\<\phi|d_k(\omega)\>
     \<d_k(\omega)|\psi\>d\mu(\omega)
\]
we obtain a splitting
\begin{equation}\label{eq:splitting}
  E_X=E^+_X-E^-_X
\end{equation}
of an arbitrary symmetric sesquilinear form measure into a
difference of two positive sesquilinear form measures. Despite the
above notation, this splitting is not canonical, and a different
choice of the operator $D$ typically yields a different
decomposition. (The choice of $\mu$ is less important: it only affects the normalization of the vectors $d_k(\omega)$.) However, by Corollary~\ref{cor:Tpm}, \emph{given}
the choice of $D$, there is a unique splitting with the stated
properties. In particular, the $\bddlin^1(H_D)$-valued extension $\tilde{E}$ (cf.\ Remark \ref{rem:HD}) has a canonical splitting into $\bddlin^1_+(H_D)$-valued operator measures. Also, if $E$ is already positive in the beginning, then the process used in the proof of the decomposition only gives $T^+=T$ and
$E^+=E$.

Let then $E:\,\Sigma\to S(V)$ be an arbitrary SFM.
\begin{definition}
The family $(E^{(k)})_{k=0}^3$ of positive SFMs
$E^{(k)}:\,\Sigma\to PS(V)$ is a {\em decomposition of $E$ (into
positive parts)} if
$$
E=\sum_{k=0}^3 i^k E^{(k)}.
$$
\end{definition}
From eqs.\ \eqref{eq:saParts} and \eqref{eq:splitting} one sees
easily that for any SFM $E$ there exists a decomposition of $E$
into positive parts.

\section{Dilations}\label{sec:dilations}

\begin{definition}\label{def:spectral dilation} Let $K$ be a Hilbert
 space, $F:\Sigma\to \mathcal L(K)$
a spectral measure, and $W\in\mathcal L(K)$ a unitary operator whose
spectrum $\sigma(W)$ is contained in $\{1,-1,i,-i\}$. Let $J:V\to
K$ be a linear map. We say that the quadruple $(K,F,W,J)$ is a
{\em (spectral $W$-)dilation} of a SFM $E:\Sigma\to S(V)$ if the
following conditions hold:
\begin{enumerate}
\item $\<J\phi|F(X)WJ\psi\>=E_X(\phi,\psi)$ for all $X\in\Sigma$ and $\phi,\,\psi\in V$,
\item $WF(X)=F(X)W$ for all $X\in\Sigma$.
\item the linear span of the set $\{W^k F(X) J\phi\,|\,k=0,1,2,3,\,X\in\Sigma,\,\phi\in V\}$ is dense in $K$.
\end{enumerate}
\end{definition}

For $k\in\{0,1,2,3\}$, let $K_k$ be the eigenspace of $W$
corresponding to $i^k$ (define $K_k=\{0\}$ if
$i^k\notin\sigma(W)$), $I_k$ the identity of $K_k$, $P_k$ the
projection of $K$ onto $K_k$, $J_k:=P_k\circ J$,
$F_k:\Sigma\to\bddlin(K_k)$ the restriction $F_k(X):=F(X)|_{K_k}$,
and $E^{(k)}:\,\Sigma\to S(V)$ the positive SFM defined by
\begin{equation}\label{eq:EXk}
  E_X^{(k)}(\phi,\psi):=\<J_k\phi|F_k(X)J_k\psi\>.
\end{equation}

\begin{theorem}\label{thm:dil}
Let $E$ be a SFM. Any dilation $(K,F,W,J)$ of $E$ defines by
\eqref{eq:EXk} a decomposition $(E^{(k)})_{k=0}^3$ of $E$ into
positive parts. Conversely, for any decomposition
$(E^{(k)})_{k=0}^3$, there exists a dilation $(K,F,W,J)$ such that
\eqref{eq:EXk} holds.
\end{theorem}

In particular, any SFM has a spectral $W$-dilation. In the
situation of Theorem~\ref{thm:dil}, we say that
$(E^{(k)})_{k=0}^3$ is the {\em decomposition of $E$ associated to
the dilation $(K,F,W,J)$}.

\begin{proof}
Given a dilation $(K,F,W,J)$, it follows from \ref{def:spectral
dilation}(2) that each $K_k$ is invariant under $F(X)$, and
$(K_k,F_k,I_k,J_k)$ is a spectral dilation of $E^{(k)}$. Then
\ref{def:spectral dilation}(1) implies that $(E^{(k)})_{k=0}^3$ is
a decomposition of $E$.

Conversely, let $(E^{(k)})_{k=0}^3$ be a decomposition of $E$.
Then each $E^{(k)}$ is a positive SFM, for which there exists a
spectral dilation of the form $(K_k,F_k,I_k,J_k)$ by Theorem~3.6
of \cite{HPY}. Define $K:=K_0\oplus K_1\oplus K_2\oplus K_3$,
$F(X):=F_0(X)\oplus F_1(X)\oplus F_2(X)\oplus F_3(X)$,
$W:=I_0\oplus(i I_1)\oplus(-I_2)\oplus(-iI_3)$, and $J:=J_0\oplus
J_1\oplus J_2\oplus J_3$. To check that $(K,F,W,J)$ is a dilation
of $E$, conditions \ref{def:spectral dilation}(1) and
\ref{def:spectral dilation}(2) are clear and \ref{def:spectral
dilation}(3) follows from
$\lin\{W^l\,|\,l=0,\ldots,3\}=\lin\{P_k\,|\,k=0,\ldots,3,\}$. It
is also clear that \eqref{eq:EXk} holds.
\end{proof}

Let $\mathcal M=(K,F,W,J)$ and  $\mathcal M'=(K',F',W',J')$ be two
dilations of $E$. The quantities $K_k'$, $F_k'$, $P_k'$, $I_k'$
and $J_k'$ related to $\mathcal M'$ are defined as before in the
obvious way.

\begin{definition}
The dilations $\mathcal M$ and  $\mathcal M'$ of $E$ are {\em
unitarily equivalent} if there exists a unitary map $U:\,K\to K'$
such that $UF(X)J\phi=F'(X)J'\phi$ for all $X\in\Sigma$, $\phi\in
V$ and $UW=W'U$; in particular, $UJ\phi=J'\phi$ for all $\phi\in
V$.
\end{definition}

\begin{theorem}
Two dilations $\mathcal M$ and  $\mathcal M'$ of $E$ are unitarily
equivalent if and only if the decompositions of $E$ associated to
$\mathcal M$ and $\mathcal M'$ are the same, that is,
$$
\  E^{(k)}_X(\phi,\psi)
  =\<J_k\phi|F(X) J_k\psi\>=\<J_k'\phi|F'(X) J_k'\psi\>
$$
for all $k\in\{0,1,2,3\}$, $X\in\Sigma$ and $\phi,\,\psi\in V$.
Moreover, then the relevant $U$ is unique and $UF(X)=F'(X)U$ for
all $X\in\Sigma$.
\end{theorem}

\begin{proof}
Assume first that the decompositions of $E$ associated to
$\mathcal M$ and $\mathcal M'$ are the same. Since both
$(K_k,F_k,I_k,J_k)$ and $(K_k',F_k',I_k',J_k')$ are dilations of
the positive SFM $E^{(k)}$, it follows from Theorem~3.6 of
\cite{HPY} that there is a unique unitary map $U_k:\,K_k\to K'_k$
such that $U_kF_k(X)J_k\phi=F'_k(X)J'_k\phi$ for all $X\in\Sigma$,
$\phi\in V$, and $U_kF_k(X)=F'_k(X)U_k$ for all $X\in\Sigma$. Then
$U:=\sum_{k=0}^3 U_k P_k$ has the desired properties.

Suppose conversely that $\mathcal M$ and $\mathcal M'$ are
unitarily equivalent. Since $UF(X)W^kF(Y)J\phi=UW^kF(X\cap
Y)J\phi=W'^kUF(X\cap Y)J\phi=W'^kF'(X\cap
Y)J'\phi=F'(X)W'^kF'(Y)UJ\phi=F'(X)UW^kF(Y)J\phi$, it follows from
\ref{def:spectral dilation}(3) that $UF(X)=F'(X)U$. As $UW=W'U$
implying $UP_k=P'_kU$, one sees that
\begin{eqnarray*}
\<P'_k J'\phi|F'(X) P'_k J'\psi\>&=& \<P'_k U J\phi|F'(X) P'_k U
J\psi\>=
\<U P_k  J\phi|F'(X) U P_k J\psi\> \\
&=&\<P_k J\phi|F(X) P_k J\psi\>,
\end{eqnarray*}
i.e., the associated decompositions coincide. Since
$UW^kF(X)J\phi=F'(X)W'^k J'\phi$, the uniqueness of $U$ is clear.
\end{proof}

\begin{remark}
Since $\| J\phi\|_K^2=\sum_{k=0}^3 \|J_k\phi\|_{K_k}^2 =
\sum_{k=0}^3 E^{(k)}_\Omega(\phi,\phi)$, we see that $J:V\to K$ is
injective if and only if
\begin{equation}\label{eq:strict}
  \sum_{k=0}^3 E^{(k)}_\Omega(\phi,\phi)>0\qquad
  \text{for all }\phi\in V\setminus\{0\}.
\end{equation}
This situation can always be achieved by writing $E=(E+\epsilon
E_0)-\epsilon E_0$, where $\epsilon>0$ and
$E_0:\,\Sigma\to\mathcal L(H)$ is a semispectral measure, which
automatically satisfies \eqref{eq:strict}.
\end{remark}

\begin{remark}
In analogy with the case of positive SFMs treated in \cite{HPY},
it is possible to describe a concrete representation of the
dilation $(K,F,W,J)$ associated with any decomposition
$(E^{(k)})_{k=0}^3$ of a SFM $E$ into positive parts.

Let $L^2(\Omega,\mu;\ell^2(\Z_+^4))\eqsim
L^2(\Omega,\mu;\ell^2(\Z_+)^4)$ be the usual Bochner space of
$\ell^2(\Z_+^4)$-valued functions $f=(f^{(0)},\ldots,f^{(3)})$,
where $f^{(k)}=(f^{(k)}_j)_{j=1}^{\infty}\in
L^2(\Omega,\mu;\ell^2(\Z_+))$. Given a measurable
$\mathbf{n}(\cdot)=(n_0(\cdot),\ldots,n_3(\cdot)):\Omega\to(\N\cup\{\infty\})^4$,
we denote by $L^2_{\mathbf{n}(\cdot)}(\Omega,\mu;\ell^2(\Z_+^4))$
the closed subspace consisting of the functions $f$ such that for
a.e.\ $\omega\in\Omega$, all $j$ and $k$, there holds
$f^{(k)}_j(\omega)=0$ if $j>n_k(\omega)$. This is analogous to the
``direct integral'' Hilbert space of a measurable family of
$\ell^2$ spaces of variable dimension considered in Section~5
of~\cite{HPY}; extending the notation used there we could write
\[
  L^2_{\mathbf{n}(\cdot)}(\Omega,\mu;\ell^2(\Z_+^4))=\int_{\Omega}^{\oplus}(\ell^2)^4_{{\bf
n}(\omega)}d\mu(\omega).
\]
Let then $E=\sum_{k=0}^3 i^k E^{(k)}$ be a SFM. By the
construction of Section~\ref{sec:diag} (or Theorem~4.5
of~\cite{HPY}), the positive SFMs $E^{(k)}$ have representations
\[
  E^{(k)}_X(\phi,\psi)
  =\int_X\sum_{j\in\Z_+}\<\phi|d^{(k)}_j(\omega)\>\<d^{(k)}_j(\omega)|\psi\>d\mu(\omega),
\]
where $d^{(k)}_j(\omega)=D^{-1}g^{(k)}_j(\omega)$, and the
$g^{(k)}_j(\omega)$ are as the $g_j(\omega)$ in
Corollary~\ref{cor:diag}. We now fix a specific
$\mathbf{n}(\cdot)$ by setting
$n_k(\omega):=\sup\{j\in\Z_+:d^{(k)}_j(\omega)\neq 0\}$ (with
$\sup\emptyset:=0$), and define
\[\begin{split}
  & K :=L^2_{\mathbf{n}(\cdot)}(\Omega,\mu;\ell^2(\Z_+^4)),\qquad
  F(X)f:=1_X f,\qquad Wf:=(f^{(0)},if^{(1)},-f^{(2)},-if^{(3)}),\\
  &(J\phi)(\omega):=
  (\<d_j^{(0)}(\omega)|\phi\>,\ldots,\<d_{j}^{(3)}(\omega)|\phi\>)_{j=1}^{\infty}.
\end{split}\]
The conditions \ref{def:spectral dilation}(1) and
\ref{def:spectral dilation}(2) of a dilation follow from simple
algebra. The density requirement \ref{def:spectral dilation}(3) is
a consequence of the fact that the component dilations
$(K_k,F_k,I_k,J_k)$, $k=0,\ldots,3$, are dilations of the positive
parts $E^{(k)}$ of $E$ by Theorem~5.1 of \cite{HPY}.
\end{remark}

\end{document}